\documentclass[11pt]{article}

\usepackage{comment}
\usepackage{color}
\usepackage{amssymb}
\usepackage{amsmath}
\usepackage[english]{babel}
\usepackage{fancyhdr}

\newtheorem{thm}{Theorem}
\newtheorem{lemma}{Lemma}

\newtheorem{cor}{Corollary}

\newtheorem{example}{Example}

\def\ep{\hfill $\Diamond$}
\def\proof{{\it Proof.}\ }
\def\J{{\mathfrak J}\, }

\def\D{{\mathfrak D}\, }

\begin{document}
\begin{center}
{\Large\bf An approach between the multiplicative and additive
structure of a Jordan ring }

\vspace{.2in}
{\bf Bruno L. M. Ferreira}
\vspace{.2in}

Technological Federal University of Paran\'{a},\\
Professora Laura Pacheco Bastos Avenue, 800, 85053-510, Guarapuava, Brazil
brunoferreira@utfpr.edu.br

\vspace{.2in}

\end{center}

{\bf keywords:} $n-$multiplicative maps, $n-$multiplicative derivations,
additivity, Jordan rings.


\begin{abstract}
Let $\J$ and $\J^{'}$ be Jordan rings. We prove under some conditions that if $\J$ contains a nontrivial idempotent, then $n-$multiplicative maps and $n-$multiplicative derivations from $\J$ to $\J^{'}$ are additive maps.

\end{abstract}

\section{Introduction}
In what follows, we use $(x, y, z) = (xy)z - x(yz)$ and $[x, y] = xy-yx$ to denote the
associator of elements $x, y, z$ and the commutator of elements $x, y$ in a not necessarily
associative ring $\J$.

For those readers who are not familiar with this language we recommend \cite{shest}
``Rings that are nearly associative, by K. A. Zhevlakov, A. M. Slin'ko, I. P. Shestakov
and A. I. Shirshov, translated by Harry F. Smith, Academic Press, New York, 1982".

Let $X = \left\{x_i\right\}_{i\in \mathbb{N}}$ be an arbitrary set of variables. A \textit{nonassociative monomial of degree $1$} is any element of $X$. Given a natural number $n > 1$, a \textit{nonassociative monomial of degree $n$} is an expression of the form $(u)(v)$, where $u$ is a nonassociative monomial of some degree $i$ and $v$ a nonassociative monomial of degree $n-i$.


Let $\J$ and $\J'$ be two rings and $\varphi: \J \longrightarrow \J'$ a bijective map of $\J$
onto $\J'$. We call $\varphi$ an \textit{$n-$multiplicative map} of $\J$ onto $\J'$ if for all
nonassociative monomials $m = m(x_1, \cdots , x_n)$ of degree $n$
$$\varphi(m(x_1, \cdots , x_n)) = m(\varphi(x_1), \cdots , \varphi(x_n))$$
for all $x_1, \ldots, x_n \in \J$.
If $\varphi(xy) = \varphi(x)\varphi(y)$ for all $x, y \in \J$, we just
say that $\varphi$ is a \textit{multiplicative map}. And if $\varphi(xyx) = \varphi(x)\varphi(y)\varphi(x)$ for all
$x, y \in \J$, then we call $\varphi$ a \textit{Jordan semi-triple multiplicative map}. 

Similarly, a map $d : \J \longrightarrow \J$ is called a \textit{$n-$multiplicative derivation} of $\J$ if
$$d(m(x_1, \cdots , x_n)) = \sum_{i=1}^{n} m(x_1,\cdots, d(x_i), \cdots , x_n)$$
for all nonassociative monomials $m = m(x_1, \ldots, x_n)$ of degree $n$ and arbitrary elements $x_1, \ldots , x_n \in \J$. 

If $d(xy) = d(x)y + xd(y)$ for all $x, y \in \J$, we just
say that $d$ is a \textit{multiplicative derivation} of $\J$. And if $d(xyx) = d(x)(yx) + xd(y)x + (xy)d(x)$ for all
$x, y \in \J$, then we call $d$ a \textit{Jordan triple multiplicative derivation}. 

A ring $\J$ is said to be \textit{Jordan} if $(x^2, y, x) = 0$ and $[x,y] = 0$ for all $x, y \in \J$. 
A Jordan ring $\J$ is called \textit{$k-$torsion free} if $k x = 0$ implies $x = 0$, for
any $x \in \J$, where $k$ is a positive integer.

A nonzero element $e \in \J$ is called an \textit{idempotent} if $e^2=ee = e$ and a
\textit{nontrivial idempotent} if it is a nonzero idempotent and different from the multiplicative
identity element of $\J$. 

Let us consider a Jordan ring $\J$ with a nontrivial idempotent $e$.

Let $\J = \J_{1} \oplus \J_{\frac{1}{2}} \oplus \J_{0}$ be the Peirce decomposition of $\J$ with respect to $e$,  where $\J_{i} = \left\{x_i ~ | ~ ex_i = ix_i\right\}$, $i = 0, \frac{1}{2}, 1$, satisfying the following multiplicative relations:

\begin{flushleft}
\ \ \ \ \ $\J_0 \J_0 \subseteq \J_0$; \ \ \  $\J_1 \J_1 \subseteq \J_1$; \ \ \  $\J_1 \J_0 = 0$; \ \ \ 
$(\J_1 \oplus \J_0)\J_\frac{1}{2} \subseteq \J_\frac{1}{2}$; \ \ \ 
$\J_\frac{1}{2} \J_\frac{1}{2} \subseteq \J_1 \oplus \J_0.$
\end{flushleft}

In studying preservers on algebras or rings, one usually assumes additivity in advance. Recently, however, a growing number of papers began investigating
preservers that are not necessarily additive. Characterizing the interrelation between the multiplicative and additive structures of a ring or algebra in an interesting topic. The first result about the additivity of maps on rings was given by Martindale III \cite{Mart}. He established a condition on a ring $\mathfrak{R}$ such that every multiplicative isomorphism on $\mathfrak{R}$ is additive. Ferreira and Ferreira \cite{Ferr} also considered this question in the context of $n-$multiplicative maps on alternative rings satisfying Martindale's conditions. Ferreira and Nascimento proved the additivity of multiplicative derivations \cite{bruth}. 

In the present paper we consider the similar Ferreira and Ferreira's problems in the context of Jordan rings. We investigate the problem of when a $n$-multiplicative isomorphism and a $n$-multiplicative derivation is additive for the class of Jordan rings.


\section{$n$-multiplicative isomorphism}
It will be convenient for us to change notation at this point.
Henceforth the ring $\J$ will be $2-$torsion free and the nonassociative monomial $m$ of degree $n$ will be an expression of the form
$$m(x_1, x_2, \ldots , x_{n-1},  x_n) := x_1(x_2(\cdots(x_{n-1}x_n)\cdots)).$$
When the first $i$ variables in the non-associative monomial $m$ assume equal values, we will denote by
$$m(\underbrace{z,z,\ldots,z}_{i},x_{i+1},\ldots,x_n) := \xi_z(x_{i+1}, \ldots, x_n).$$


The main technique we will use is the following argument which will be termed a ``standard argument". 
Suppose, $x, y, s \in \J$ are such that $\varphi(s) = \varphi(x) + \varphi(y)$. Multiplying this equality by $\varphi(t_i)$, $(i = 1, 2, \ldots , n-1)$, we get
\begin{eqnarray*}
\varphi(t_1)(\varphi(t_2)(\cdots(\varphi(t_{n-1})\varphi(s))\cdots)) &=& \varphi(t_1)(\varphi(t_2)(\cdots(\varphi(t_{n-1})\varphi(x))\cdots)) \\&+& \varphi(t_1)(\varphi(t_2)(\cdots(\varphi(t_{n-1})\varphi(y))\cdots))
\end{eqnarray*}
so
\begin{eqnarray*}
m(\varphi(t_1),\varphi(t_2),\ldots,\varphi(t_{n-1}),\varphi(s)) &=& m(\varphi(t_1),\varphi(t_2),\ldots,\varphi(t_{n-1}),\varphi(x))\\&+& m(\varphi(t_1),\varphi(t_2),\ldots,\varphi(t_{n-1}),\varphi(y)) 
\end{eqnarray*}

It follows that
\begin{eqnarray*}
\varphi \big(m(t_1, t_2, \ldots , t_{n-1}, s)\big) = \varphi \big(m(t_1, t_2, \ldots , t_{n-1}, x)\big) + \varphi \big(m(t_1, t_2, \ldots , t_{n-1}, y)\big).
\end{eqnarray*}

Moreover, if
\begin{eqnarray*}
\varphi \big(m(t_1, t_2, \ldots , t_{n-1}, x)\big) + \varphi \big(m(t_1, t_2, \ldots , t_{n-1}, y)\big) &=& \varphi \big(m(t_1, t_2, \ldots , t_{n-1}, x) \\&+& m(t_1, t_2, \ldots , t_{n-1}, y)\big)
\end{eqnarray*}

then by injectivity of $\varphi$, we have that

$$m(t_1, t_2, \ldots , t_{n-1}, s) = m(t_1, t_2, \ldots , t_{n-1}, x) + m(t_1, t_2, \ldots , t_{n-1}, y).$$

The main result of this section reads as follows.

\begin{thm}\label{t11} Let $\J$ and $\J'$ be Jordan rings and $e$ a non-trivial idempotent in $\J$. Let $\J = \J_1 \oplus \J_{\frac{1}{
2}}\oplus J_0$ be the Peirce decomposition of $\J$ with respect to $e$. If $\J$ satisfies the following conditions:
\begin{enumerate}
\item[\it(i)] Let $a_i \in \J_i(i = 1, 0)$. If $t_{\frac{1}{2}}a_i= 0$ for all $t_{\frac{1}{2}} \in \J_{\frac{1}{2}}$, then $a_i = 0$;
\item[\it(ii)] Let $a_0 \in \J_0$. If $t_0a_0  = 0$ for all $t_0 \in \J_0$, then $a_0 = 0$;
\item[\it(iii)] Let $a_{\frac{1}{2}}\in \J_{\frac{1}{2}}$. If $t_0a_{\frac{1}{2}} = 0$ for all $t_0 \in \J_0$, then $a_{\frac{1}{2}}= 0$.
\end{enumerate}
Then every $n-$multiplicative isomorphism from $\J$ onto $\J'$ is additive.
\end{thm}

The proof is organized in a serie of Lemmas.

\begin{lemma}\label{lema1}
$\varphi(0) = 0.$
\end{lemma}
\proof Since $\varphi$ is surjective, there exists $x \in \J$ such that $\varphi(x) = 0$. Therefore $\varphi(\xi_0(x)) = \xi_{\varphi(0)}(\varphi(x)) = \xi_{\varphi(0)}(0) = 0.$ \ep

\begin{lemma}\label{lema2}
Let $a_i \in \J_i$, $i = 1, \frac{1}{2}, 0$. Then $\varphi(a_1 + a_{\frac{1}{2}} + a_0) = \varphi(a_1) + \varphi(a_{\frac{1}{2}}) + \varphi(a_0).$
\end{lemma}
\proof Since $\varphi$ is surjective, we can find an element $s = s_1 + s_{\frac{1}{2}} + s_0 \in \J$ such that 
\begin{eqnarray}\label{l01}
\varphi(s) = \varphi(a_1) + \varphi(a_{\frac{1}{2}})+ \varphi(a_0).
\end{eqnarray}
For $e$, applying a standard argument to (\ref{l01}), we get 
\begin{eqnarray*}
\varphi \big(\xi_{2e}(s)\big) &=& \varphi \big(\xi_{2e}(a_1)\big) + \varphi \big(\xi_{2e}(a_{\frac{1}{2}})\big) + \varphi \big(\xi_{2e}(a_0)\big) =  \varphi(2^{n-1}a_1) + \varphi(a_{\frac{1}{2}}) + \varphi(0) \\&=&  \varphi(2^{n-1}a_1) + \varphi(a_{\frac{1}{2}}). 
\end{eqnarray*}
Since $\varphi(\xi_{2e}(s)) = \varphi(2^{n-1}s_1 + s_{\frac{1}{2}})$, we have 
\begin{eqnarray}\label{l02}
\varphi(2^{n-1}s_1 + s_{\frac{1}{2}}) = \varphi(2^{n-1}a_1) + \varphi(a_{\frac{1}{2}}).
\end{eqnarray}
Now for $t_0 \in \J_0$ applying the standard argument to (\ref{l02}), we have
\begin{eqnarray*}
\varphi \big(\xi_{2e}(t_0,2^{n-1}s_1 + s_{\frac{1}{2}})\big)= \varphi \big(\xi_{2e}(t_0,a_1)\big) + \varphi \big(\xi_{2e}(t_0,a_{\frac{1}{2}})\big) =\varphi(t_0a_{\frac{1}{2}}).
\end{eqnarray*}
Hence $$\varphi(t_0s_{\frac{1}{2}}) = \varphi(t_0a_{\frac{1}{2}}).$$
Therefore, $t_0s_{\frac{1}{2}} = t_0a_{\frac{1}{2}}$ for every $t_0 \in \J_0$. It follows from the item $(iii)$ of Theorem \ref{t11}  that $s_{\frac{1}{2}} = a_{\frac{1}{2}}.$
With a similar argument we can show that $s_1 = a_1$.

Now only remains show that $s_0 = a_0$. For this let $t_0 \in \J_0$, applying the standard argument to (\ref{l01}), we get
\begin{eqnarray}\label{l04}
\nonumber \varphi \big(\xi_{t_0}(s)\big) &=& \varphi \big(\xi_{t_0}(a_1)\big) + \varphi \big(\xi_{t_0}(a_{\frac{1}{2}})\big)\nonumber + \varphi \big(\xi_{t_0}(a_0)\big) \nonumber = \varphi \big( 0 \big) + \varphi \nonumber\big(\xi_{t_0}(a_{\frac{1}{2}})\big) + \varphi \big(\xi_{t_0}(a_0)\big) 
\\&=& \varphi \big(\xi_{t_0}(a_{\frac{1}{2}})\big) + \varphi \big(\xi_{t_0}(a_0)\big). 
\end{eqnarray}
For $t_{\frac{1}{2}} \in \J_{\frac{1}{2}}$, applying the standard argument to (\ref{l04}), we have
\begin{eqnarray}\label{l05}
\nonumber \varphi \big(\xi_{2e}(t_{\frac{1}{2}},\xi_{t_0}(s))\big) &=& \varphi \big(\xi_{2e}(t_{\frac{1}{2}},\xi_{t_0}(a_{\frac{1}{2}}))\big) \nonumber + \varphi \big(\xi_{2e}(t_{\frac{1}{2}},\xi_{t_0}(a_0)\big) \nonumber = \varphi \big(2^{n-2}(t_{\frac{1}{2}}\xi_{t_0}(a_{\frac{1}{2}}))_1\big)\nonumber \\&+& \varphi \big(t_{\frac{1}{2}}\xi_{t_0}(a_{0})\big),
\end{eqnarray}
where $(t_{\frac{1}{2}}\xi_{t_0}(a_{\frac{1}{2}}))_1 \in \J_{1}$ and $\displaystyle t_{\frac{1}{2}}\xi_{t_0}(a_{0}) \in \J_{\frac{1}{2}}$. Now, for $t^{'}_0 \in \J_0$, applying the standard argument to (\ref{l05}),
we have that 
\begin{eqnarray*}
\varphi \big(\xi_{2e}(t^{'}_0, \xi_e(t_{\frac{1}{2}},\xi_{t_0}(s)))\big) &=& \varphi \big(\xi_{2e}(t^{'}_0,2^{n-2}(t_{\frac{1}{2}}\xi_{t_0}(a_{\frac{1}{2}}))_1)\big) + \varphi \big(\xi_{2e}(t^{'}_0,t_{\frac{1}{2}}\xi_{t_0}(a_{0}))\big) \\&=&  \varphi(0) + \varphi \big(\xi_{2e}(t^{'}_0,t_{\frac{1}{2}}\xi_{t_0}(a_{0}))\big) \\&=& \varphi \big(\xi_{2e}(t^{'}_0,t_{\frac{1}{2}}\xi_{t_0}(a_{0}))\big)\\&=&  \varphi \big(t^{'}_0(t_{\frac{1}{2}}\xi_{t_0}(a_{0}))\big) .
\end{eqnarray*}
Since $\xi_{2e}(t^{'}_0, \xi_e(t_{\frac{1}{2}},\xi_{t_0}(s))) = \displaystyle t^{'}_0(t_{\frac{1}{2}}\xi_{t_0}(s_{0})),$ we have
$$t^{'}_0(t_{\frac{1}{2}}\xi_{t_0}(s_{0})) = t^{'}_0(t_{\frac{1}{2}}\xi_{t_0}(a_{0})).$$
It follows from the items $(i)$, $(ii)$ and $(iii)$ of Theorem \ref{t11} that $s_0 = a_0$. Thus $s = a_1 + a_{\frac{1}{2}} + a_0$. \ep

\begin{lemma}\label{lema3}
Let $a_{\frac{1}{2}}, b_{\frac{1}{2}} \in \J_{\frac{1}{2}}$ and $a_0 \in \J_0$. Then $\varphi(a_{\frac{1}{2}}a_0 + b_{\frac{1}{2}}) = \varphi(a_{\frac{1}{2}}a_0) +\varphi(b_{\frac{1}{2}})$. 
\end{lemma}
\proof We note
$$\xi_{2e}(2e+a_{\frac{1}{2}}, a_0 + b_{\frac{1}{2}}) = \xi_{2e}(2e, a_{\frac{1}{2}}a_0 + b_{\frac{1}{2}}) + \xi_{2e}(b_{\frac{1}{2}}, a_{\frac{1}{2}}).$$
Using Lemma \ref{lema2} we have $\varphi(\xi_{2e}(2e+a_{\frac{1}{2}}, a_0 + b_{\frac{1}{2}})) = \varphi(\xi_{2e}(2e, a_{\frac{1}{2}}a_0 + b_{\frac{1}{2}})) + \varphi(\xi_{2e}(b_{\frac{1}{2}}, a_{\frac{1}{2}}))$ because $\xi_{2e}(2e, a_{\frac{1}{2}}a_0 + b_{\frac{1}{2}}) \in \J_{\frac{1}{2}}$ and $\xi_{2e}(b_{\frac{1}{2}}, a_{\frac{1}{2}}) \in \J_1$.
Consequently, by Lemmas \ref{lema1} and \ref{lema2}, we have
\begin{eqnarray*}
\varphi(\xi_{2e}(2e, a_{\frac{1}{2}}a_0 + b_{\frac{1}{2}})) &+& \varphi(\xi_{2e}(a_{\frac{1}{2}},b_{\frac{1}{2}})) = \varphi(\xi_{2e}(2e+a_{\frac{1}{2}}, a_0 + b_{\frac{1}{2}}))\\&&= \xi_{\varphi(2e)}(\varphi(2e+a_{\frac{1}{2}}), \varphi(a_0 + b_{\frac{1}{2}})) \\&&= \xi_{\varphi(2e)}(\varphi(2e)+\varphi(a_{\frac{1}{2}}), \varphi(a_0) + \varphi(b_{\frac{1}{2}}))\\&&= \xi_{\varphi(2e)}(\varphi(2e)+\varphi(a_{0}))+ \xi_{\varphi(2e)}(\varphi(2e) , \varphi(b_{\frac{1}{2}})) \\&&+ \xi_{\varphi(2e)}(\varphi(a_{\frac{1}{2}}) , \varphi(a_0)) + \xi_{\varphi(2e)}(\varphi(a_{\frac{1}{2}}) , \varphi(b_{\frac{1}{2}}))\\&&= \varphi(0) + \varphi(\xi_{2e}(2e , b_{\frac{1}{2}})) + \varphi(\xi_{2e}(a_{\frac{1}{2}} , a_0)) \\&&+ \varphi(\xi_{2e}(a_{\frac{1}{2}} , b_{\frac{1}{2}})).
\end{eqnarray*}
Thus, $\varphi(\xi_{2e}(2e, a_{\frac{1}{2}}a_0 + b_{\frac{1}{2}})) = \varphi(\xi_{2e}(2e , b_{\frac{1}{2}})) + \varphi(\xi_{2e}(a_{\frac{1}{2}} , a_0))$, that is $\varphi(a_{\frac{1}{2}}a_0 + b_{\frac{1}{2}}) = \varphi(a_{\frac{1}{2}}a_0) +\varphi(b_{\frac{1}{2}}).$ \ep
\begin{lemma}\label{lema4}
Let $a_{\frac{1}{2}}, b_{\frac{1}{2}} \in \J_{\frac{1}{2}}$. Then $\varphi(a_{\frac{1}{2}}+ b_{\frac{1}{2}})= \varphi(a_{\frac{1}{2}}) + \varphi(b_{\frac{1}{2}})$.
\end{lemma}
\proof Let $s = s_1 + s_{\frac{1}{2}} + s_0 \in \J$ such that 
\begin{eqnarray}\label{idmeio}
\varphi(s) = \varphi(a_{\frac{1}{2}}) + \varphi(b_{\frac{1}{2}}).
\end{eqnarray}
For $t_0 \in \J_0$, applying the standard argument to (\ref{idmeio}) and using Lemma \ref{lema3}, we have
$$\varphi(\xi_{t_0}(s,t_0)) = \varphi(\xi_{t_0}(a_{\frac{1}{2}},t_0)) + \varphi(\xi_{t_0}(b_{\frac{1}{2}},t_0)) = \varphi(\xi_{t_0}(a_{\frac{1}{2}},t_0) + \xi_{t_0}(b_{\frac{1}{2}},t_0)).$$
Hence $\xi_{t_0}(s,t_0) = \xi_{t_0}(a_{\frac{1}{2}},t_0) + \xi_{t_0}(b_{\frac{1}{2}},t_0).$ By $(ii)$ and $(iii)$ of Theorem \ref{t11}, we obtain that $s_0 = 0$ and $s_{\frac{1}{2}} = a_{\frac{1}{2}} + b_{\frac{1}{2}}.$
Now for $t_{\frac{1}{2}} \in \J_{\frac{1}{2}}$ and $e \in \J_1$, applying the standard argument to (\ref{idmeio}) again, we have
\begin{eqnarray}\label{idmeio2}
\varphi(\xi_{2e}(s_1,t_{\frac{1}{2}}) + \xi_{2e}(s_{\frac{1}{2}}, t_{\frac{1}{2}})) = \varphi(\xi_{2e}(s, t_{\frac{1}{2}})) = \varphi(\xi_{2e}(a_{\frac{1}{2}}, t_{\frac{1}{2}})) + \varphi(\xi_{2e}(b_{\frac{1}{2}}, t_{\frac{1}{2}})).
\end{eqnarray}
For $u_0$, applying the standard argument to (\ref{idmeio2}), we get that
\begin{eqnarray*}
\varphi(\xi_{u_0}(\xi_{2e}(s_1,t_{\frac{1}{2}}))) = \varphi(\xi_{u_0}(\xi_{2e}(a_{\frac{1}{2}}, t_{\frac{1}{2}}))) + \varphi(\xi_{u_0}(\xi_{2e}(b_{\frac{1}{2}}, t_{\frac{1}{2}}))) = \varphi(0) + \varphi(0) = 0.
\end{eqnarray*}
Hence $\xi_{u_0}(\xi_{2e}(s_1,t_{\frac{1}{2}})) = 0$ for every $u_0 \in \J_0$ and $t_{\frac{1}{2}} \in \J_{\frac{1}{2}}$, it follows from the items $(iii)$ and $(i)$ of Theorem \ref{t11} that $s_1 = 0$. 
Thus, $s= s_{\frac{1}{2}} = a_{\frac{1}{2}} + b_{\frac{1}{2}}.$ \ep
\begin{lemma}\label{lema5}
Let $a_{1}, b_{1} \in \J_{1}$. Then $\varphi(a_{1}+ b_{1})= \varphi(a_{1}) + \varphi(b_{1})$.
\end{lemma}
\proof Let $s = s_1 + s_{\frac{1}{2}} + s_0 \in \J$ such that 
\begin{eqnarray}\label{idum}
\varphi(s) = \varphi(a_1) + \varphi(b_1).
\end{eqnarray}
For $t_0 \in \J_0$, applying the standard argument to (\ref{idum}), we have
$$\varphi(\xi_{t_0}(s_{\frac{1}{2}},t_0) + \xi_{t_0}(s_0,t_0))=\varphi(\xi_{t_0}(s,t_0)) = \varphi(\xi_{t_0}(a_{1},t_0)) + \varphi(\xi_{t_0}(b_{1},t_0)) = \varphi(0) + \varphi(0) = 0.$$
Thus, $\xi_{t_0}(s_{\frac{1}{2}},t_0) + \xi_{t_0}(s_0,t_0) = 0$ for every $t_0 \in J_0$. Since $\xi_{t_0}(s_{\frac{1}{2}},t_0) \in \J_{\frac{1}{2}}$ and $\xi_{t_0}(s_0,t_0) \in \J_0$, we get that $\xi_{t_0}(s_{\frac{1}{2}},t_0) = 0$ and $\xi_{t_0}(s_0,t_0) = 0$. By $(ii)$ and $(iii)$ of Theorem \ref{t11}, we obtain that $s_{\frac{1}{2}} = 0$ and $s_0 = 0$.
Now for $t_{\frac{1}{2}} \in \J_{\frac{1}{2}}$ and $e \in \J_1$, applying the standard argument to (\ref{idum}) again and using Lemma \ref{lema4}, we have
\begin{eqnarray}\label{idum2}
\nonumber \varphi(\xi_{2e}(s_1,t_{\frac{1}{2}})) = \varphi(\xi_{2e}(s, t_{\frac{1}{2}})) &=& \varphi(\xi_{2e}(a_{1}, t_{\frac{1}{2}})) + \varphi(\xi_{2e}(b_{1}, t_{\frac{1}{2}})) \\&=& \varphi(\xi_{2e}(a_{1}, t_{\frac{1}{2}}) + \xi_{2e}(b_{1}, t_{\frac{1}{2}})).
\end{eqnarray}
Hence $\xi_{2e}(s_1,t_{\frac{1}{2}}) = \xi_{2e}(a_{1}, t_{\frac{1}{2}}) + \xi_{2e}(b_{1}, t_{\frac{1}{2}})$ for every $t_{\frac{1}{2}} \in \J_{\frac{1}{2}}$, it follows from the item $(i)$ of Theorem \ref{t11} that $s_1 = a_1 + b_1$. 
Thus, $s= s_{1} = a_{1} + b_{1}.$ \ep
\begin{lemma}\label{lema6}
Let $a_{0}, b_{0} \in \J_{0}$. Then $\varphi(a_{0}+ b_{0})= \varphi(a_{0}) + \varphi(b_{0})$.
\end{lemma}
\proof Let $s = s_1 + s_{\frac{1}{2}} + s_0 \in \J$ such that 
\begin{eqnarray}\label{idzero}
\varphi(s) = \varphi(a_0) + \varphi(b_0).
\end{eqnarray}
For $e \in \J_1$, applying the standard argument to (\ref{idzero}), we have
$$\varphi(\xi_{2e}(s_{\frac{1}{2}},2e) + \xi_{2e}(s_1,2e))=\varphi(\xi_{2e}(s,2e)) = \varphi(\xi_{2e}(a_{0},2e)) + \varphi(\xi_{2e}(b_{0},2e)) = \varphi(0) + \varphi(0) = 0.$$
Therefore, $s_{\frac{1}{2}} + 2^{n-1} s_1 = 0$. Hence $s_{\frac{1}{2}}=0$ and $s_1 = 0$. 
Now for $t_{\frac{1}{2}} \in \J_{\frac{1}{2}}$ and $e \in \J_1$, applying the standard argument to (\ref{idzero}) again and using Lemma \ref{lema4}, we have
\begin{eqnarray}\label{idzero2}
\nonumber \varphi(\xi_{2e}(s_0,t_{\frac{1}{2}})) = \varphi(\xi_{2e}(s, t_{\frac{1}{2}})) &=& \varphi(\xi_{2e}(a_{0}, t_{\frac{1}{2}})) + \varphi(\xi_{2e}(b_{0}, t_{\frac{1}{2}})) \\&=& \varphi(\xi_{2e}(a_{0}, t_{\frac{1}{2}}) + \xi_{2e}(b_{0}, t_{\frac{1}{2}})).
\end{eqnarray}
Hence $\xi_{2e}(s_0,t_{\frac{1}{2}}) = \xi_{2e}(a_{0}, t_{\frac{1}{2}}) + \xi_{2e}(b_{0}, t_{\frac{1}{2}})$ for every $t_{\frac{1}{2}} \in \J_{\frac{1}{2}}$, it follows from the item $(i)$ of Theorem \ref{t11} that $s_0 = a_0 + b_0$. 
Thus, $s= s_{0} = a_{0} + b_{0}.$ \ep
\vspace{0.5cm}

Now we are ready to prove our first result.
\vspace{0.5cm}

\noindent{\bf Proof of Theorem \ref{t11}.}

\noindent \proof Let $a= a_1 + a_{\frac{1}{2}} + a_0$, $b= b_1 + b_{\frac{1}{2}} + b_0$. By Lemmas \ref{lema2}, \ref{lema4}, \ref{lema5} and \ref{lema6} we have
\begin{eqnarray*}
\varphi(a+b) &=& \varphi((a_1 + b_1)+(a_{\frac{1}{2}} + b_{\frac{1}{2}})+(a_0 + b_0))
\\&=& \varphi(a_1 + b_1)+\varphi(a_{\frac{1}{2}} + b_{\frac{1}{2}})+\varphi(a_0 + b_0)
\\&=& \varphi(a_1) + \varphi(b_1)+\varphi(a_{\frac{1}{2}}) + \varphi(b_{\frac{1}{2}})+\varphi(a_0) + \varphi(b_0)
\\&=& \varphi(a_1 + a_{\frac{1}{2}} + a_0) + \varphi(b_1 + b_{\frac{1}{2}} + b_0)
\\&=& \varphi(a) + \varphi(b).
\end{eqnarray*}
That is, $\varphi$ is additive on $\J$. \ep

\section{$n$-multiplicative derivation}

We now investigate the problem of when a $n$-multiplicative derivations is additive for the class of Jordan rings. 

For this purpose we will assume that the Jordan ring $\J$ is $\left\{2,(n-1)\right\}$-torsion free for $n\geq2$ where $n$ is degree of the nonassociative monomial $m = m(x_1, \cdots , x_n)$.

Let $d$ be a $n$-multiplicative derivation of Jordan ring $\J$. If we put $d(e) = a_1 +a_{\frac{1}{2}} + a_0$, then $\displaystyle d(m(e,e,\cdots,e)) = \sum_{i=1}^{n} m(e,\cdots, d(e), \cdots , e)  = n a_1  + a_{\frac{1}{2}}.$ 
Since $d(m(e,e,\cdots,e)) = d(e)$ then $(n-1)a_1 - a_0 = 0.$ Thus, $a_1 = a_0=0$ and $d(e)= a_{\frac{1}{2}}$. By [77 page of \cite{Sch}] we have 
$$\mathcal{D}_{y,z}(x) = [L_{y},L_{z}] + [L_y, R_z] + [R_y,R_z]$$
is a derivation for all $y, z \in \J$. In particular, if $y = a_{\frac{1}{2}}$ and $z = 4e$, then $\mathcal{D}_{y,z}(e) = 3d(e)$. In deed,
\begin{eqnarray*}
\mathcal{D}_{a_{\frac{1}{2}},4e}(e) &=& ([L_{a_{\frac{1}{2}}},L_{4e}] + [L_{a_{\frac{1}{2}}}, R_{4e}] + [R_{a_{\frac{1}{2}}},R_{4e}])(e)
\\&=& L_{a_{\frac{1}{2}}}L_{4e}(e) - L_{4e}L_{a_{\frac{1}{2}}}(e) + L_{a_{\frac{1}{2}}}R_{4e}(e) - R_{4e}L_{a_{\frac{1}{2}}}(e) \\&+& R_{a_{\frac{1}{2}}}R_{4e}(e) - R_{4e}R_{a_{\frac{1}{2}}}(e) \\&=& a_{\frac{1}{2}}(4ee) - 4e(a_{\frac{1}{2}}e) + a_{\frac{1}{2}}(e4e) - (a_{\frac{1}{2}}e)4e + (e4e)a_{\frac{1}{2}} - (ea_{\frac{1}{2}})4e \\&=& 2a_{\frac{1}{2}} - a_{\frac{1}{2}} + 2a_{\frac{1}{2}} - a_{\frac{1}{2}} + 2a_{\frac{1}{2}} - a_{\frac{1}{2}} \\&=& 3a_{\frac{1}{2}}, 
\end{eqnarray*}
so $(\mathcal{D}_{y,z}-3d)(e) = 0$.
Let $\D = \mathcal{D}_{y,z}-3d$, note that $\D$ is additive if and only if $d$ is additive, since $\mathcal{D}_{y,z}$ is additive.
Furthermore, observe that $\D$ is a multiplicative derivation such that $\D(e) = 0$.
\vspace{0.5cm}

The next is the main result of this section. Its proof shares the same
outline as that of Theorem \ref{t11} but it needs different technique.

\begin{thm}\label{t22} Let $\J$ be a Jordan ring with a non-trivial idempotent $e$. Let $\J = \J_1 \oplus \J_{\frac{1}{
2}}\oplus J_0$ be the Peirce decomposition of $\J$ with respect to $e$. If $\J$ satisfies the conditions of Theorem \ref{t11},
then every $n-$multiplicative derivation $d$ from $\J$ is additive.
\end{thm}

The proof will be organized in a series of auxiliary lemmas.

\begin{lemma}\label{d1}
$\D(0) = 0$. 
\end{lemma}
\proof Note that $\D(0) = \D(\xi_0(0)) = 0$ 
\ep 

\begin{lemma}\label{d2}
$\D(\J_i) \subseteq \J_i$ for $i = 1, \frac{1}{2}, 0.$
\end{lemma}
\proof Let $a_1 \in \J_1$ as $\D(e) = 0$ we have $\D(a_1) = \D(\xi_e(a_1)) = \xi_e(\D(a_1)).$
If we express $\D(a_1) = \D(a_1)_1 + \D(a_1)_{\frac{1}{2}} + \D(a_1)_0$ it follows that 
$$\D(a_1) = \D(a_1)_1 + \D(a_1)_{\frac{1}{2}} + \D(a_1)_0 = \xi_e(\D(a_1)) = \D(a_1) + \frac{1}{2^{n-1}} \D(a_1)_{\frac{1}{2}}.$$
Thus, $\D(a_1) = 0 = \D(a_1)_{\frac{1}{2}}$ and $\D(a_1) = \D(a_1)_1$.

Let $a_0 \in \J_0$ as $\D(e) = 0$ we have $0 = \D(\xi_e(a_0)) = \xi_e(\D(a_0)) = \D(a_0)_1 + \frac{1}{2^{n-1}}\D(a_0)_{\frac{1}{2}}.$
If we express $\D(a_0) = \D(a_0)_1 + \D(a_0)_{\frac{1}{2}} + \D(a_0)_0$ it follows that $\D(a_0)_1 = 0 = \D(a_0)_{\frac{1}{2}}.$
Thus, $\D(a_0) = \D(a_0)_0$.

Let $a_{\frac{1}{2}} \in \J_{\frac{1}{2}}$ we have
\begin{eqnarray*}
\D(a_{\frac{1}{2}}) &=& \D(\xi_{2e}(a_{\frac{1}{2}})) \\&=& \underbrace{\xi_{\D(2e)}(2e,\ldots,2e,a_{\frac{1}{2}}) + \xi_{2e}(\D(2e), 2e, \ldots , 2e,a_{\frac{1}{2}}) + \cdots + \xi_{2e}(2e,\ldots,\D(2e),a_{\frac{1}{2}})}_{x_{\frac{1}{2}} \in \J_{\frac{1}{2}}} \\&+& \xi_{2e}(2e,\ldots,2e,\D(a_{\frac{1}{2}})). 
\end{eqnarray*}

If we express $\D(a_{\frac{1}{2}}) = \D(a_{\frac{1}{2}})_1 + \D(a_{\frac{1}{2}})_{\frac{1}{2}} + \D(a_{\frac{1}{2}})_0$ it follows that $(2^{n-1} - 1)\D(a_{\frac{1}{2}})_1 = 0 = \D(a_{\frac{1}{2}})_{0}$ and $x_{\frac{1}{2}} = 0$ once that $\xi_{2e}(2e,\ldots,2e,\D(a_{\frac{1}{2}})) = 2^{n-1}\D(a_{\frac{1}{2}})_1 + \D(a_{\frac{1}{2}})_{\frac{1}{2}}$.
Thus, $\D(a_{\frac{1}{2}}) = \D(a_{\frac{1}{2}})_{\frac{1}{2}}$.

\begin{lemma}\label{d3}
Let $a_1 \in \J_1$ and $a_{\frac{1}{2}} \in \J_{\frac{1}{2}}$. Then $\D(2^{n-1}a_{\frac{1}{2}}) = 2^{n-1}\D(a_{\frac{1}{2}})$, $\D(2e) = 0$ and $\D(2^{n-1}a_1) = 2^{n-1}\D(a_1).$ 
\end{lemma}
\proof Let $a_{\frac{1}{2}} \in \J_{\frac{1}{2}}$ as $\D(e) = 0$ and $\D(2^{n-1}a_{\frac{1}{2}}) \in \J_{\frac{1}{2}}$ then
$$\D(a_{\frac{1}{2}}) = \D(\xi_e(2^{n-1}a_{\frac{1}{2}})) = \xi_e(\D(2^{n-1}a_{\frac{1}{2}})) = \frac{1}{2^{n-1}}\D(2^{n-1}a_{\frac{1}{2}}).$$
Thus, $\D(2^{n-1}a_{\frac{1}{2}}) = 2^{n-1}\D(a_{\frac{1}{2}})$.

Now note that $(n-1)\D(2e)a_{\frac{1}{2}} = \D(\xi_{2e}(a_{\frac{1}{2}})) - \xi_{2e}(\D(a_{\frac{1}{2}})) = \D(a_{\frac{1}{2}}) - \D(a_{\frac{1}{2}}) = 0.$
By item $(i)$ of Theorem \ref{t22} we have $\D(2e) = 0.$

Let $a_{1} \in \J_{1}$ as $\D(2e) = 0$ and $\D(a_{1}) \in \J_{1}$ then
$$\D(2^{n-1}a_1) = \D(\xi_{2e}(a_1)) = \xi_{2e}(\D(a_1)) = 2^{n-1}\D(a_1).$$ \ep

\begin{lemma}\label{d4}
Let $a_1 \in \J_1$, $a_{\frac{1}{2}} \in \J_{\frac{1}{2}}$ and $a_0 \in \J_0$. Then $\D(a_1 + a_{\frac{1}{2}} + a_0 ) = \D(a_1) + \D(a_{\frac{1}{2}}) + \D(a_0).$
\end{lemma}
\proof Consider $\D(a_1 + a_{\frac{1}{2}} + a_0 ) = d_1 + d_{\frac{1}{2}} + d_0$, by Lemma \ref{d3} we get
\begin{eqnarray}\label{ostres}
\nonumber \D(2^{n-1}a_1 + a_{\frac{1}{2}}) &=& \D(\xi_{2e}(a_1 + a_{\frac{1}{2}} + a_0 )) = \xi_{2e}(\D(a_1 + a_{\frac{1}{2}} + a_0 ))\\&=& \xi_{2e}(d_1 + d_{\frac{1}{2}} + d_0) = 2^{n-1}d_1 + d_{\frac{1}{2}}.
\end{eqnarray}
Let $t_0 \in \J_0$ then 
\begin{eqnarray}\label{dif}
\nonumber t_0d_{\frac{1}{2}} = \xi_{2e}(t_0,\D(2^{n-1}a_1 + a_{\frac{1}{2}})) &=& \D(\xi_{2e}(t_0, 2^{n-1}a_1 + a_{\frac{1}{2}})) - \xi_{2e}(\D(t_0),2^{n-1}a_1 + a_{\frac{1}{2}}) \\&=& \D(t_0a_{\frac{1}{2}}) - \D(t_0)a_{\frac{1}{2}} = t_0\D(a_{\frac{1}{2}}). 
\end{eqnarray}
Thus, by item (iii) of Theorem \ref{t22}  we have $\D(a_{\frac{1}{2}}) = d_{\frac{1}{2}}.$
Let $t_{\frac{1}{2}} \in \J_{\frac{1}{2}}$ by (\ref{ostres}) we have
\begin{eqnarray}\label{maisum}
\nonumber \xi_{2e}(t_{\frac{1}{2}}, 2^{n-1}d_1) + \xi_{2e}(t_{\frac{1}{2}}, d_{\frac{1}{2}}) &=& \xi_{2e}(t_{\frac{1}{2}}, 2^{n-1}d_1 + d_{\frac{1}{2}}) = \xi_{2e}(t_{\frac{1}{2}}, \D(2^{n-1}a_1+a_{\frac{1}{2}})) \nonumber \\&=& \D(\xi_{2e}(t_{\frac{1}{2}}, 2^{n-1}a_1 + a_{\frac{1}{2}})) - \xi_{2e}(\D(t_{\frac{1}{2}}), 2^{n-1}a_1 + a_{\frac{1}{2}})\nonumber \\&=& \nonumber \D(\xi_{2e}(t_{\frac{1}{2}}, 2^{n-1}a_1) + \xi_{2e}(t_{\frac{1}{2}}, a_{\frac{1}{2}})) \\&-& \xi_{2e}(\D(t_{\frac{1}{2}}), 2^{n-1}a_1) - \xi_{2e}(\D(t_{\frac{1}{2}}), a_{\frac{1}{2}})
\end{eqnarray} 
As $\xi_{2e}(t_{\frac{1}{2}}, d_{\frac{1}{2}}) , \xi_{2e}(\D(t_{\frac{1}{2}}), a_{\frac{1}{2}}) \in \J_1 \oplus \J_0$ follows that for $u_0 \in \J_0$ by (\ref{maisum}) we get
\begin{eqnarray*}
u_0\xi_{2e}(t_{\frac{1}{2}}, 2^{n-1}d_1) &=& \xi_{2e}(u_0, \xi_{2e}(t_{\frac{1}{2}} , 2^{n-1}d_1)) = \xi_{2e}(u_0, \xi_{2e}(t_{\frac{1}{2}} , 2^{n-1}d_1) + \xi_{2e}(t_{\frac{1}{2}}, d_{\frac{1}{2}})) \\&=& \xi_{2e}(u_0, \D(\xi_{2e}(t_{\frac{1}{2}},2^{n-1}a_1) + \xi_{2e}(t_{\frac{1}{2}}, a_{\frac{1}{2}}))) - \xi_{2e}(u_0,\xi_{2e}(\D(t_{\frac{1}{2}}),\\&&2^{n-1}a_1)) \\&=& u_0\xi_{2e}(t_{\frac{1}{2}} , \D(2^{n-1}a_1)).
\end{eqnarray*}
Therefore by item $(iii)$ of Theorem \ref{t22} and Lemma \ref{d3} we have $d_1 = \D(a_1).$ 

Finally, we show that $d_0 = \D(a_0)$. Let $e \in \J_1$, $t_{\frac{1}{2}} \in \J_{\frac{1}{2}}$ and $h_0, t_0 \in \J_0$. We have
$$\xi_{h_0}(h_0, 2e, t_{\frac{1}{2}}, t_0, \D(a_1 + a_{\frac{1}{2}} + a_0)) = \xi_{h_0}(h_0, 2e, t_{\frac{1}{2}}, t_0, d_1 + d_{\frac{1}{2}} + d_0) = \xi_{h_0}(h_0, t_{\frac{1}{2}}, t_0, d_0).$$

On the other hand using the identity of $n-$multiplicative derivation we get $$\xi_{h_0}(h_0, 2e, t_{\frac{1}{2}}, t_0, \D(a_1 + a_{\frac{1}{2}}+ a_0)) = \xi_{h_0}(h_0, t_{\frac{1}{2}}, t_0, \D(a_0)).$$
Thus $\xi_{h_0}(h_0, t_{\frac{1}{2}}, t_0, d_0)=\xi_{h_0}(h_0, t_{\frac{1}{2}}, t_0, \D(a_0))$. Now by items $(iii), (i)$ and $(ii)$ of Theorem \ref{t22} we have $d_0 = \D(a_0).$

Therefore $\D(a_1 + a_{\frac{1}{2}} + a_0) = \D(a_1) + \D(a_{\frac{1}{2}}) + \D(a_0).$
 \ep
\begin{lemma}\label{d5}
Let $a_{\frac{1}{2}}, b_{\frac{1}{2}} \in \J_{\frac{1}{2}}$. Then $\D(a_{\frac{1}{2}} + b_{\frac{1}{2}} ) = \D(a_{\frac{1}{2}}) + \D(b_{\frac{1}{2}}).$
\end{lemma}
\proof Note that $\xi_{2e}((2e + a_{\frac{1}{2}}),(2e+b_{\frac{1}{2}})) = \xi_{2e}(a_{\frac{1}{2}} + b_{\frac{1}{2}} +2^2e + a_{\frac{1}{2}}b_{\frac{1}{2}}).$
Since $a_{\frac{1}{2}} + b_{\frac{1}{2}} = \xi_{2e}(a_{\frac{1}{2}} + b_{\frac{1}{2}}) \in \J_{\frac{1}{2}}$ and $\xi_{2e}(2^2e + a_{\frac{1}{2}}b_{\frac{1}{2}}) \in \J_{1}$ by Lemmas \ref{d3} and \ref{d4} we get
\begin{eqnarray*}
\D(\xi_{2e}(a_{\frac{1}{2}} + b_{\frac{1}{2}})) + \D(\xi_{2e}(2^2e + a_{\frac{1}{2}}b_{\frac{1}{2}}))&=& \xi_{2e}(\D(a_{\frac{1}{2}} + b_{\frac{1}{2}})) + \xi_{2e}(\D(2^2e + a_{\frac{1}{2}}b_{\frac{1}{2}}))\\&=& \xi_{2e}(\D(a_{\frac{1}{2}} + b_{\frac{1}{2}}) + \D(2^2e + a_{\frac{1}{2}}b_{\frac{1}{2}}))\\&=& \xi_{2e}(\D(a_{\frac{1}{2}} + b_{\frac{1}{2}} + 2^2e + a_{\frac{1}{2}}b_{\frac{1}{2}}))\\&=& \D(\xi_{2e}(a_{\frac{1}{2}} + b_{\frac{1}{2}} + 2^2e + a_{\frac{1}{2}}b_{\frac{1}{2}}))\\&=& \D(\xi_{2e}((2e + a_{\frac{1}{2}}),(2e + b_{\frac{1}{2}})))\\&=& \xi_{2e}(\D(2e + a_{\frac{1}{2}}),(2e + b_{\frac{1}{2}}))) \\&+& \xi_{2e}((2e + a_{\frac{1}{2}}),\D(2e + b_{\frac{1}{2}})))\\&=&\xi_{2e}(\D(2e) + \D(a_{\frac{1}{2}}),(2e + b_{\frac{1}{2}})) \\&+& \xi_{2e}((2e + a_{\frac{1}{2}}),\D(2e) + \D(b_{\frac{1}{2}}))\\&=& \xi_{2e}(\D(a_{\frac{1}{2}}),2e) + \xi_{2e}(\D(a_{\frac{1}{2}}), b_{\frac{1}{2}}) \\&+& \xi_{2e}(2e , \D(b_{\frac{1}{2}})) + \xi_{2e}(a_{\frac{1}{2}}, \D(b_{\frac{1}{2}})) \\&=& \D(\xi_{2e}(a_{\frac{1}{2}}, b_{\frac{1}{2}})) + \xi_{2e}(\D(a_{\frac{1}{2}}), 2e) \\&+& \xi_{2e}(2e, \D(b_{\frac{1}{2}})).
\end{eqnarray*}
Observe that $\D(\xi_{2e}(2^2e + a_{\frac{1}{2}}b_{\frac{1}{2}})), \D(\xi_{2e}(a_{\frac{1}{2}}, b_{\frac{1}{2}})) \in \J_1$ and
 $\D(\xi_{2e}(a_{\frac{1}{2}} + b_{\frac{1}{2}})),\xi_{2e}(\D(a_{\frac{1}{2}}), 2e) + \xi_{2e}(2e, \D(b_{\frac{1}{2}})) \in \J_{\frac{1}{2}}$ it follows that $$\D(a_{\frac{1}{2}} + b_{\frac{1}{2}} ) = \D(\xi_{2e}(a_{\frac{1}{2}} + b_{\frac{1}{2}})) = \xi_{2e}(\D(a_{\frac{1}{2}}), 2e) + \xi_{2e}(2e, \D(b_{\frac{1}{2}})) = \D(a_{\frac{1}{2}}) + \D(b_{\frac{1}{2}}).$$ \ep

\begin{lemma}\label{d6}
Let $a_{i}, b_{i} \in \J_{i}$, $(i = 1, 0)$. Then $\D(a_{i} + b_{i} ) = \D(a_{i}) + \D(b_{i}).$
\end{lemma}
\proof Let $a_i, b_i \in \J_i$ $(i = 1, 0)$ and $t_{\frac{1}{2}} \in \J_{\frac{1}{2}}$, by Lemma \ref{d5} we have
\begin{eqnarray*}
\D(a_i+b_i)t_{\frac{1}{2}}=\xi_{2e}(\D(a_i+b_i), t_{\frac{1}{2}}) &=& \D(\xi_{2e}((a_i+b_i),t_{\frac{1}{2}})) - \xi_{2e}((a_i+b_i),\D(t_{\frac{1}{2}}))\\&=&
\D(\xi_{2e}(a_i,t_{\frac{1}{2}})) + \D(\xi_{2e}(b_i,t_{\frac{1}{2}}))\\&-&\xi_{2e}(a_i,\D(t_{\frac{1}{2}})) - \xi_{2e}(b_i,\D(t_{\frac{1}{2}}))\\&=& \xi_{2e}(\D(a_i),t_{\frac{1}{2}}) + \xi_{2e}(\D(b_i),t_{\frac{1}{2}})\\&=& (\D(a_i)+\D(b_i))t_{\frac{1}{2}}.
\end{eqnarray*}
Therefore, by item $(i)$ of Theorem \ref{t22} we get $\D(a_{i} + b_{i} ) = \D(a_{i}) + \D(b_{i})$. \ep
\vspace{0.5cm}

Now we are in a position to show that $\D$ preserves addition.
\vspace{0.5cm}

\noindent{\bf Proof of Theorem \ref{t22}.}

\noindent \proof Let $a= a_1 + a_{\frac{1}{2}} + a_0$, $b= b_1 + b_{\frac{1}{2}} + b_0$. By Lemmas \ref{d4}-\ref{d6} we have
\begin{eqnarray*}
\D(a+b) &=& \D((a_1 + b_1)+(a_{\frac{1}{2}} + b_{\frac{1}{2}})+(a_0 + b_0))
\\&=& \D(a_1 + b_1)+\D(a_{\frac{1}{2}} + b_{\frac{1}{2}})+\D(a_0 + b_0)
\\&=& \D(a_1) + \D(b_1)+\D(a_{\frac{1}{2}}) + \D(b_{\frac{1}{2}})+\D(a_0) + \D(b_0)
\\&=& \D(a_1 + a_{\frac{1}{2}} + a_0) + \D(b_1 + b_{\frac{1}{2}} + b_0)
\\&=& \D(a) + \D(b).
\end{eqnarray*}
That is, $\D$ is additive on $\J$. \ep

\vspace{0.2in}

The following two examples show that there are non-trivial noncommutative Jordan algebra and Jordan algebra respectively that satisfy the conditions of the Theorem \ref{t11}.

\begin{example}\label{ex1} Let $\mathfrak{F}$ be a field of characteristic different from $2$, $\mathfrak{J}$ a four dimensional algebra over $\mathfrak{F}$ and a basis $\{e_{11}, e_{10}, e_{01}, e_{00}\}$ with the multiplication table given by: $e_{ij}e_{kl}=\delta _{jk} e_{il}~(i,j,k,l=1,2)$, where $\delta _{jk}$ is the Kronecker delta. We can verify that $\mathfrak{J}$ is a Jordan algebra. In fact, $\mathfrak{J}$ is an associative algebra where $e_{11}$ and $e_{00}$ are orthogonal idempotents such that $e=e_{11} + e_{00}$ is the unity element of $\mathfrak{J}$. Moreover, if $\mathfrak{J}=\mathfrak{J}_{1}\oplus
\mathfrak{J}_{\frac{1}{2}}\oplus \mathfrak{J}_{0}$ is the Peirce decomposition of $\mathfrak{J}$, relative to $e_{11}$, then we have $\mathfrak{J}_{1}=\mathfrak{F}e_{11}, \mathfrak{J}_{\frac{1}{2}} = \mathfrak{F}e_{10} + \mathfrak{F}e_{01}, \mathfrak{J}_{0}=\mathfrak{F}e_{00}$. From a direct calculation, we can verify that $\mathfrak{J}$
satisfies the conditions {\it (i)-(iii)} of the Theorem \ref{t11}.
\end{example}

\begin{example} Let $\mathfrak{K}$ be the algebra obtained from the associative algebra $\mathfrak{J}$, in Example \ref{ex1}, on replacing the product $xy$ by $x\cdot y=\frac{1}{2}(xy+yx)$. We can verify that $\mathfrak{K}$ is a Jordan algebra where $e_{11}$ and $e_{00}$ are orthogonal idempotents such that $e=e_{11}+e_{00}$ is the unity element of $\mathfrak{K}$. Moreover, if $\mathfrak{K}=\mathfrak{K}_{1}\oplus \mathfrak{K}_{\frac{1}{2}}\oplus \mathfrak{K}_{0}$ is the Peirce decomposition of $\mathfrak{K}$, relative to $e_{11}$, then we have $\mathfrak{K}_{i}=\mathfrak{F}e_{ii}~(i=1,0)$ and $\mathfrak{K}_{\frac{1}{2}}=\mathfrak{F}e_{10}+\mathfrak{F}e_{01}$. From a direct calculation, we can verify that the algebra $\mathfrak{K}$ satisfies the conditions {\it (i)-(iii)} of the Theorem \ref{t11}.
\end{example}

\section{Corollaries}
In this section, we will give some consequences of our main results.

\begin{cor}
Let $\J$ satisfy the conditions of Theorem \ref{t11}. Then any multiplicative map $\varphi$ of $\J$
onto an arbitrary Jordan ring $\J'$ is additive.
\end{cor}


In the case of unital Jordan rings we have

\begin{cor}
Let $\J$ be a unital Jordan ring and $e$ a non-trivial idempotent in $\J$. Let $\J'$ be a Jordan
ring. Let $\J = \J_1 \oplus \J_{\frac{1}{2}}\oplus \J_0$ be the Peirce decomposition of $\J$ with respect to $e$. If $\J$ satisfies the
following conditions:
\begin{enumerate}
\item[\it(i)] Let $a_i \in \J_i(i = 1, 0)$. If $t_{\frac{1}{2}}a_i= 0$ for all $t_{\frac{1}{2}} \in \J_{\frac{1}{2}}$, then $a_i = 0$,
\end{enumerate}
then every multiplicative map from $\J$ onto $\J'$ is additive.
\end{cor}

\begin{cor}
Let $\J$ be a unital Jordan ring and $e$ a non-trivial idempotent in $\J$. Let $\J = \J_1 \oplus \J_{\frac{1}{2}}\oplus \J_0$ be the Peirce decomposition of $\J$ with respect to $e$. If $\J$ satisfies the
following condition:
\begin{enumerate}
\item[\it(i)] Let $a_i \in \J_i(i = 1, 0)$. If $t_{\frac{1}{2}}a_i= 0$ for all $t_{\frac{1}{2}} \in \J_{\frac{1}{2}}$, then $a_i = 0$,
\end{enumerate}
then every multiplicative derivation from $\J$ is additive.
\end{cor}

\end{document}